\documentclass[11pt, a4paper, oneside]{amsart}

\RequirePackage{amsmath}
\RequirePackage{bm}
\RequirePackage{amssymb}
\RequirePackage{upref}
\RequirePackage{amsthm}
\RequirePackage{enumerate}
\RequirePackage{pb-diagram}
\RequirePackage{amsfonts}
\RequirePackage[mathscr]{eucal}
\RequirePackage{verbatim}
\RequirePackage{xr}
\RequirePackage{graphicx}
\usepackage{calc}
\usepackage{xspace}

\newcommand{\e}{\epsilon}

\newcommand{\f}{\varphi}
\newcommand{\h}{\eta}

\newcommand{\m}{\mu}
\newcommand{\n}{\nu}

\newcommand{\s}{\sigma}
\newcommand{\x}{\xi}

\newcommand{\C}{\varGamma}

\newcommand{\F}{\varPhi}

\newcommand{\di}[1]{#1\nobreakdash-\hspace{0pt}dimensional}
\newcommand{\nbdd}{\nobreakdash--}
\newcommand{\nbd}{\nobreakdash-\hspace{0pt}}

\newcommand{\ces}[1]{$C^#1$\nbd{estimates}}

\newcommand{\fm}[1]{F_{|_{M_#1}}}
\newcommand{\fmo}[1]{F_{|_{#1}}}
\newcommand{\fu}[3]{#1\hspace{0pt}_{|_{#2_#3}}}
\newcommand{\fv}[2]{#1\hspace{0pt}_{|_{#2}}}

\newcommand{\so}{{\mc S_0}}

\newcommand{\const}{\tup{const}}

\newcommand{\msp[1]}[1]{\mspace{#1mu}}
\newcommand{\low}[1]{{\hbox{}_{#1}}}

\newcommand{\R}[1][n+1]{{\protect\mathbb R}^{#1}}

\newcommand{\N}{{\protect\mathbb N}}

\newcommand{\eR}{\stackrel{\lower1ex \hbox{\rule{6.5pt}{0.5pt}}}{\msp[3]\R[]}}
\newcommand{\eN}{\stackrel{\lower1ex \hbox{\rule{6.5pt}{0.5pt}}}{\msp[1]\N}}
\newcommand{\eO}{\stackrel{\lower1ex
\hbox{\rule{6pt}{0.5pt}}}{\msc O}}

\DeclareMathOperator{\graph}{graph}

\newcommand\pa{\partial}
\newcommand\pde[2]{\frac {\partial#1}{\partial#2}}
\newcommand\pd[3]{\frac {\partial#1}{\partial#2^#3}}   
\newcommand\pdc[3]{\frac {\partial#1}{\partial#2_#3}}   
\newcommand\pdm[4]{\frac {\partial#1}{\partial#2_#3^#4}}   
\newcommand\pddc[4]{\frac {{\partial\hskip0.15em}^2#1}{\partial {#2_
#3}\,\partial{#2_#4}}} 

\newcommand\sql[1][u]{\sqrt{1-|D#1|^2}}
\newcommand{\un}{\infty}
\newcommand{\A}{\forall}

\newcommand{\set}[2]{\{\,#1\colon #2\,\}}
\newcommand{\uu}{\cup}
\newcommand{\ii}{\cap}
\newcommand{\uuu}{\bigcup}

\newcommand{\uud}{ \stackrel{\lower 1ex \hbox {.}}{\uu}}
\newcommand{\uuud}[1]{ \stackrel{\lower 1ex \hbox {.}}{\uuu_{#1}}}
\newcommand\su{\subset}

\newcommand{\sminus}[1][28]{\raise 0.#1ex\hbox{$\scriptstyle\setminus$}}

\newcommand{\abs}[1]{\lvert#1\rvert}

\newcommand{\norm}[1]{\lVert#1\rVert}

\newcommand{\nnorm}[1]{| \mspace{-2mu} |\mspace{-2mu}|#1| \mspace{-2mu}
|\mspace{-2mu}|}
\newcommand{\spd}[2]{\protect\langle #1,#2\protect\rangle}

\newcommand\ch[3]{\varGamma_{#1#2}^#3}
\newcommand\cha[3]{{\bar\varGamma}_{#1#2}^#3}

\newcommand{\riem}[4]{R_{#1#2#3#4}}
\newcommand{\riema}[4]{{\bar R}_{#1#2#3#4}}

\newcommand{\tbf}{\textbf}
\newcommand{\tit}{\textit}

\newcommand{\tup}{\textup}

\newcommand{\mc}{\protect\mathcal}
\newcommand{\msc}{\protect\mathscr}

\providecommand{\bysame}{\makebox[3em]{\hrulefill}\thinspace}

\newcommand{\ci}{\cite}
\newcommand{\bib}{\bibitem}

\newcommand{\bt}{\begin{thm}}
\newcommand{\bl}{\begin{lem}}
\newcommand{\bc}{\begin{cor}}
\newcommand{\bd}{\begin{definition}}
\newcommand{\bpp}{\begin{prop}}
\newcommand{\br}{\begin{rem}}
\newcommand{\bn}{\begin{note}}
\newcommand{\be}{\begin{ex}}
\newcommand{\bes}{\begin{exs}}
\newcommand{\bb}{\begin{example}}
\newcommand{\bbs}{\begin{examples}}
\newcommand{\ba}{\begin{axiom}}

\newcommand{\et}{\end{thm}}
\newcommand{\el}{\end{lem}}
\newcommand{\ec}{\end{cor}}
\newcommand{\ed}{\end{definition}}
\newcommand{\epp}{\end{prop}}
\newcommand{\er}{\end{rem}}
\newcommand{\en}{\end{note}}
\newcommand{\ee}{\end{ex}}
\newcommand{\ees}{\end{exs}}
\newcommand{\eb}{\end{example}}
\newcommand{\ebs}{\end{examples}}
\newcommand{\ea}{\end{axiom}}

\newcommand{\bp}{\begin{proof}}
\newcommand{\ep}{\end{proof}}
\newcommand{\eps}{\renewcommand{\qed}{}\end{proof}}

\newcommand{\bal}{\begin{align}}

\newcommand{\bi}[1][1.]{\begin{enumerate}[\upshape #1]}
\newcommand{\bia}[1][(1)]{\begin{enumerate}[\upshape #1]}
\newcommand{\bin}[1][1]{\begin{enumerate}[\upshape\bfseries #1]}
\newcommand{\bir}[1][(i)]{\begin{enumerate}[\upshape #1]}
\newcommand{\bic}[1][(i)]{\begin{enumerate}[\upshape\hspace{2\cma}#1]}
\newcommand{\bis}[2][1.]{\begin{enumerate}[\upshape\hspace{#2\parindent}#1]}
\newcommand{\ei}{\end{enumerate}}

\newcommand\ndots{\raise 0.47ex \hbox {,}\hskip0.06em\cdots %
     \raise 0.47ex \hbox {,}\hskip0.06em} 

\newcommand{\q}{\quad}
\newcommand{\qq}{\qquad}

\newcommand{\hp}{\hphantom}

\newcommand\nd{\noindent}

\newskip\Csmallskipamount                                                
\Csmallskipamount=\smallskipamount
\newskip\Cmedskipamount
\Cmedskipamount=\medskipamount
\newskip\Cbigskipamount
\Cbigskipamount=\bigskipamount

\newcommand\cvs{\vspace\Csmallskipamount}   
\newcommand\cvm{\vspace\Cmedskipamount}

\newskip\csa
\csa=\smallskipamount

\newskip\cma
\cma=\medskipamount

\newskip\cba
\cba=\bigskipamount

\newdimen\spt
\newcommand\citem{\cvs\advance\itemno by
1{(\romannumeral\the\itemno})\hskip3pt}
\newcommand{\bitem}{\cvm\nd\advance\itemno by
1{\bf\the\itemno}\hspace{\cma}}

\newcount\itemno
\newcommand{\las}[1]{\label{S:#1}}

\newcommand{\lae}[1]{\label{E:#1}}
\newcommand{\lat}[1]{\label{T:#1}}
\newcommand{\lal}[1]{\label{L:#1}}
\newcommand{\lad}[1]{\label{D:#1}}

\newcommand{\lap}[1]{\label{P:#1}}
\newcommand{\lar}[1]{\label{R:#1}}

\newcommand{\rs}[1]{Section~\ref{S:#1}}

\newcommand{\rt}[1]{Theorem~\ref{T:#1}}
\newcommand{\rl}[1]{Lemma~\ref{L:#1}}
\newcommand{\rd}[1]{Definition~\ref{D:#1}}

\newcommand{\rp}[1]{Proposition~\ref{P:#1}}
\newcommand{\rr}[1]{Remark~\ref{R:#1}}

\newcommand{\re}[1]{\eqref{E:#1}}

\newskip\thmskip
\thmskip=\parindent

\newskip\hsk
\newenvironment{hinw}{\labelsep=0pt\begin{list}{}{\labelsep=0pt\itemindent=0pt\labelwidth=0pt\leftmargin=\parindent\rightmargin=0pt\partopsep=\cba}%
\item\it\nopagebreak\nopagebreak}%
{\end{list}}

\newcommand\bh{\begin{hinw}}
\newcommand{\eh}{\end{hinw}}

\newtheoremstyle{normal}
  {\cba}
  {\cba}
  {}
  {\thmskip}
  {\bfseries}
  {.}
  {\hsk}
  {}

\newtheoremstyle{abschnitt}
  {\cba}
  {\cba}
  {}
  {\thmskip}
  {\bfseries}
  {.}
  {\hsk}
  {}

\newtheoremstyle{italic}
  {\cba}
  {\cba}
  {\itshape}
  {\thmskip}
  {\bfseries}
  {.}
  {\hsk}
  {}

\newtheoremstyle{aufgaben}
  {\cba}
  {\cba}
  {}
  {}
  {\normalsize\bfseries}
  {.}
  {\hsk}
  {}

\newtheoremstyle{break}
  {\cba}
  {\cba}
  {\itshape}
  {}
  {\bfseries}
  {.}
  {\newline}
  {}

\swapnumbers
\theoremstyle{italic}
\newtheorem{thm}[subsection]{Theorem}
\newtheorem{lem}[subsection]{Lemma}
\newtheorem{prop}[subsection]{Proposition}
\newtheorem{cor}[subsection]{Corollary}

\theoremstyle{normal}
\newtheorem{rem}[subsection]{Remark}
\newtheorem{definition}[subsection]{Definition}
\newtheorem{example}[subsection]{Example}
\newtheorem{examples}[subsection]{Examples}
\newtheorem{ex}[subsection]{Exercise}
\newtheorem{note}[subsection]{}
\newtheorem{axiom}[subsection]{Axiom}

\theoremstyle{aufgaben}
\newtheorem{exs}[subsection]{Exercises}

\swapnumbers

\numberwithin{equation}{section}
\numberwithin{figure}{section}

\newenvironment{textequation}[1][0.8]
{\begin{equation}
\begin{aligned}
\begin{minipage}{#1\linewidth}}
{\end{minipage}
\end{aligned}
\end{equation}
\ignorespacesafterend}

\newcommand{\btext}{\begin{textequation}}
\newcommand{\etext}{\end{textequation}}

\RequirePackage{upref}
\RequirePackage{enumerate}
\usepackage[mathscr]{eucal}

\usepackage{hyperref}
\begin{document}

\title[Hypersurfaces of prescribed  curvature]{Hypersurfaces of
prescribed  curvature in Lorentzian manifolds}

\author{Claus Gerhardt}
\address{Ruprecht-Karls-Universit\"at, Institut f\"ur Angewandte Mathematik,
Im Neuenheimer Feld 294, 69120 Heidelberg, Germany}
\email{gerhardt@math.uni-heidelberg.de}

%
\date{\today}
%


\begin{abstract} The existence of closed hypersurfaces of prescribed 
curvature  in globally hyperbolic Lorentzian manifolds is proved provided there
are barriers.
\end{abstract}
\maketitle

\tableofcontents
\setcounter{section}{-1}
\section{Introduction} 

Consider the problem of finding a closed hypersurface of prescribed curvature
$F$ in a complete \di {(n+1)} manifold $N$. To be more precise, let
$\Omega$ be a connected open subset of $N,f\in C^{2,\alpha}(\bar\Omega), F$ a smooth,
symmetric function defined in an open cone $\C\su\R[n]$, then we look for a
hypersurface
$M\su \Omega$ such that
\begin{equation}\lae{0.1}
F_{|_M}=f(x)\qq\A x\in M,
\end{equation}
where $F_{|_M}$ means that $F$ is evaluated at the vector $(\kappa_i(x))$ the
components of which are the principal curvatures of $M$. The prescribed function
$f$ should satisfy natural structural conditions, e. g. if $\varGamma$ is the
positive cone and the hypersurface $M$ is supposed to be convex, then $f$ should
be positive, but no further, merely technical, conditions should be imposed.

If $N$ is a Riemannian manifold, then the problem has been solved in the case
when $F=H$, the {\it mean curvature}, where in addition $n$ had to be small, and
$N$ conformally flat, cf. \ci{cg5}, and for curvature functions $F$ of class $(K)$,
no restrictions on $n$, cf. \ci{cg2, cg4}. We also refer to \ci{cg3}, where more
special situations are considered, and the bibliography therein.

In Lorentzian manifolds, existence results for space-like hypersurfaces of
prescribed curvature have only been proved in the case $F=H$ so far, cf. \ci{cg1},
and \ci{rb, rb2}, where the results are better than in the Riemannian case, since
no restrictions on
$n$ have to be imposed, and rather general ambient spaces can be considered.
Thus, one would hope that for curvature functions of class $(K)$ the existence
results are at least as good as in the Riemannian case, and maybe the proof a
little bit less demanding.

Unfortunately, the Lorentzian structure is only of advantage as far as the
$C^1$\nbd{estimates} are concerned, while the proof of the \ces{2} is more
difficult, if not impossible, for arbitrary functions of class $(K)$. The
complications derive from the Gau{\ss} equations, where the term stemming
from the second fundamental form of the hypersurface has the opposite sign in
the Lorentzian case as compared to the Riemannian case, which in turn leads to
an unfavourable sign in the equation for the second fundamental form used for
the a priori estimates.

We were only able to overcome these difficulties for curvature functions
belonging to a fairly large subclass of $(K)$, called $(K^*)$, which will be defined
in
\rs{1}, and which includes the Gaussian curvature.

To give a precise statement of the existence result, we need a few definitions
and assumptions. First, we assume that $N$ is a smooth, connected, \tit{globally
hyperbolic} manifold with a compact \tit{Cauchy hypersurface}, or equivalently,
that $N$ is topologically a product, $N=\R[]\times\mc{S}_0$, where $\mc S_0$ is a
compact, \di n Riemannian manifold, and there exists a Gaussian coordinate
system $(x^\alpha)_{0\le \alpha\le n}$ such that $x^0$ represents the time, the
$(x^i)_{1\le i\le n}$ are local coordinates for $\mc S_0$, where we may assume
that $\mc S_0$ is equal to the level hypersurface $\{x^0=0\}$---we don't
distinguish between $\mc S_0$ and $\{0\}\times \mc S_0$---, and such that the
Lorentzian metric takes the form
\begin{equation}\lae{0.2}
d\bar s_N^2=e^{2\psi}\{-{dx^0}^2+\sigma_{ij}(x^0,x)dx^idx^j\},
\end{equation}
where $\sigma_{ij}$ is a Riemannian metric, $\psi$ a function on $N$, and $x$ an
abbreviation for the space-like components $(x^i)$, see \ci{GR},
\ci[p.~212]{HE}, \ci[p.~252]{GRH}, and \ci[Section~6]{cg1}.

In $N$ we consider an open, connected set $\Omega$ that is bounded by two
\tit{achronal}, connected, space-like hypersurfaces $M_1\text{ and }M_2$, where
$M_1$ is supposed to lie in the \tit{past} of $M_2$.

Let $F$ be of class $(K^*)$, and $0<f\in C^{2,\alpha}(\bar\Omega)$. Then, we assume that the
boundary components $M_i$ act as barriers for $(F,f)$.

\bd\lad{0.1}
$M_2$ is an \tit{upper barrier} for $(F,f)$, if $M_2$ is strictly convex and
satisfies
\begin{equation}
\fm 2\ge f,
\end{equation}
and  $M_1$ is a \tit{lower barrier} for $(F,f)$, if at the points $\varSigma\su M_1$, where
$M_1$ is strictly convex, there holds
\begin{equation}
\fv F\varSigma\le f.
\end{equation}
$\varSigma$ may be empty.
\ed

We shall clarify in \rs 2 what \tit {convexity} means for space-like hypersurfaces.

Then, we can prove

\bt\lat{0.2}
Let $M_1$ be a lower and $M_2$ an upper barrier for $(F,f)$. Then, the problem
\begin{equation}
\fmo M= f
\end{equation}
has a strictly convex solution $M\su \bar\Omega$ of class $C^{4,\alpha}$ that can be
written as a graph over $\mc S_0$ provided there exists a strictly convex
function $\chi\in C^2(\bar\Omega)$.
\et

\br
As we shall show in \rs 2 the existence of a strictly convex function $\chi$ is
guaranteed by the assumption that the level hypersurfaces $\{x^0=\tup{const}\}$
are strictly convex in $\bar\Omega$.
\er
The paper is organized as follows: In \rs 1 we define the curvature functions of
class $(K^*)$ and examine their properties.

In \rs 2 we introduce the notations and common definitions we rely on,  state
the equations of Gau{\ss}, Codazzi, and Weingarten for space-like hypersurfaces
in \tit{pseudo-riemannian} manifolds, and analyze achronal hypersurfaces in some
detail.

In \rs 3 we look at the curvature flow associated with our problem, and the
corresponding evolution equations for the basic geometrical quantities of the
flow hypersurfaces.

In \rs 4 we prove lower order estimates for the evolution problem,  while a
priori estimates in the
$C^2$\nbd{norm} are derived in
\rs 5.

Finally, in \rs 6, we demonstrate that the evolutionary solution converges to a
stationary solution.

\section{Curvature functions}\las{1}

Let $\C_+\su\R[n]$ be the open positive cone and $F\in C^{2,\alpha}(\C_+)\ii
C^0(\bar\C_+)$ a symmetric function satisfying the condition
\begin{equation}
F_i=\pd F\kappa i>0\; ;
\end{equation}
then, $F$ can also be viewed as a function defined on the space of symmetric,
positive definite matrices $\mathscr S_+$, for, let $(h_{ij})\in \msc S_+$ with
eigenvalues $\kappa_i,\,1\le i\le n$, then define $F$ on $\msc S_+$ by
\begin{equation}
F(h_{ij})=F(\kappa_i).
\end{equation}

If we define 
\begin{align}
F^{ij}&=\pde F{h_{ij}}\\
\intertext{and}
F^{ij,kl}&=\pddc Fh{{ij}}{{kl}}
\end{align}
then, 
\begin{equation}
F^{ij}\x_i\x_j=\pdc F\kappa i \abs{\x^i}^2\q\A\, \x\in\R[n],
\end{equation}
\begin{equation}
F^{ij} \,\text{is diagonal if $h_{ij}$ is diagonal,}
\end{equation}
and
\begin{equation}\lae{1.7}
F^{ij,kl}\h_{ij}\h_{kl}=\pddc F{\kappa}ij\h_{ii}\h_{jj}+\sum_{i\ne
j}\frac{F_i-F_j}{\kappa_i-\kappa_j}(\h_{ij})^2,
\end{equation}
for any $(\h_{ij})\in \msc S$, where $\msc S$ is the space of all symmetric
matrices. The second term on the right-hand side of \re{1.7} is non-positive if
$F$ is concave, and non-negative if $F$ is convex, and has to be interpreted as a
limit if $\kappa_i=\kappa_j$.

In \ci{cg4} we defined the class $(K)$ as

\bd\lad{1.1}
A symmetric function $F\in C^{2,\alpha}(\C_+)\ii C^0(\bar\C_+)$ positively
homogeneous of degree 1 is said to be of class $(K)$ if
\begin{equation}\lae{1.8}
F_i=\pdc F{\kappa}i>0\qq \text{in}\; \C_+,
\end{equation}
\begin{equation}
F\,\text{is concave},
\end{equation}
\begin{equation}
F_{|_{\pa \C_+}}=0,
\end{equation}
and there exists a constant $c=c(F)$ such that
\begin{equation}\lae{1.11}
F^{ij,kl}\h_{ij}\h_{kl}\le cF^{-1}(F^{ij}\h_{ij})^2-F^{ik}\tilde
h^{jl}\h_{ij}\h_{kl}\qq\A\,\h\in\msc S,
\end{equation}
where $F$ is evaluated at $(h_{ij})\in\msc S_+$ and $(\tilde h^{ij})=(h_{ij})^{-1}$.
\ed

As we only  recently became aware of, inequality \re{1.11} is  valid with constant
$c=1$ if it is valid for a larger $c$.
\bl
Let $F\in C^2(\C_+)$ be a symmetric curvature function, positively
homogeneous of degree $d_0>0$ that satisfies the relations \re{1.8} and
\re{1.11}, then it fulfills \re{1.11} with constant $c=1$, i.e.
\begin{equation}\lae{1.12}
F^{ij,kl}\h_{ij}\h_{kl}\le F^{-1}(F^{ij}\h_{ij})^2-F^{ik}\tilde
h^{jl}\h_{ij}\h_{kl}\qq\A\,\h\in\msc S,
\end{equation}
or, equivalently, if we set $\hat F=\log F$,
\begin{equation}\lae{1.13}
\hat F^{ij,kl}\h_{ij}\h_{kl}\le -\hat F^{ik}\tilde
h^{jl}\h_{ij}\h_{kl}\qq\A\,\h\in\msc S.
\end{equation}
Equality holds in \re{1.12} and \re{1.13} for $(\h_{ij})=( h_{ij})$.
\el
\bp
As we have shown in \ci[Lemma 1.3 and Remark 1.4]{cg4} a symmetric curvature
function $F\in C^2(\C_+)$ satisfies inequality \re{1.11} iff
\begin{equation}\lae{1.14}
F_i\kappa_i\le F_j\kappa_j\,,\q \text{for} \q \kappa_j\le \kappa_i,
\end{equation}
and
\begin{equation}\lae{1.20}
F_{ij}\x^i\x^j\le cF^{-1}(F_i\x^i)^2-F_i\kappa_i^{-1}\abs{\x^i}^2\q \forall\,\x\in \R[n],
\end{equation}
where $F_i$, $F_{ij}$ are ordinary partial derivatives of $F$ in $\C_+$. Thus, we
have to show that \re{1.20} holds with $c=1$ for the $F$'s under consideration.

We note that $F>0$, cf. the proof of \rl{1.8} below. Let $\hat F=\log F$ and
\begin{equation}
f_{ij}=\hat F_{ij}+\hat F_i\kappa_i^{-1}\delta_{ij},
\end{equation} 
then the relation \re{1.20} is equivalent to 
\begin{equation}
f_{ij}-(c-1)\hat F_i\hat F_j\le 0.
\end{equation}
We shall demonstrate that
\begin{equation}
f_{ij}\le 0.
\end{equation}

Define $\Lambda$ by
\begin{equation}
\Lambda=\set{\lambda\in \R[]_+}{f_{ij}-\lambda\hat F_i\hat F_j\le 0},
\end{equation}
and let $\lambda_0=\inf \Lambda$.
 
$\Lambda$ is non-empty, so that the infimum is well defined and
attained. If $\lambda_0=0$, then the main part of the Lemma is proved. Thus, assume
that $\lambda_0>0$, and let $\m$ be the largest eigenvalue of
\begin{equation}
f_{ij}-\lambda_0\hat F_i\hat F_j
\end{equation}
with eigenspace $E$. Evidently, $\m$ must be zero.

Let $(\kappa^i)$ be the argument of $F$. Then, in view of the homogeneity of $F$ we
conclude
\begin{equation}\lae{1.26}
f_{ij}\kappa^j=0
\end{equation}
and
\begin{equation}
\hat F_i\kappa^i=d_0.
\end{equation}

Now, let $\h=(\h^i)\in E$, then

\begin{equation}
f_{ij}\h^j-\lambda_0\hat F_i\hat F_j\h^j=0,
\end{equation}
and, multiplying this equation with $(\kappa^i)$, we obtain
\begin{equation}
\lambda_0\hat F_i\h^i=0,
\end{equation}
i.e. $D\hat F$ is orthogonal to $E$, and
\begin{equation}
f_{ij}\h^j=0.
\end{equation}

For $0<\e<\lambda_0$ set
\begin{equation}
g_{ij}^\e=f_{ij}-(\lambda_0-\e)\hat F_i\hat F_j.
\end{equation}
Then the largest eigenvalue of $g_{ij}^\e$,  has to be positive because
of the definition of $\lambda_0$. Let $\h_\e$ be a corresponding unit eigenvector,
then, $\h_\e$ has to be orthogonal to $E$, for $E$ is also an eigenspace of
$g_{ij}^\e$; but this is impossible, since a subsequence of the $\h_\e$'s
converges to a unit vector in $E$, if $\e$ tends to zero.

Hence, we conclude that $\lambda_0=0$ and that inequality \re{1.13} is valid. Finally,
equality holds in \re{1.13} if we choose $(\h_{ij})=(h_{ij})$ in view of \re{1.26}.
\ep

Thus, it seems worth to redefine the class $(K)$.
\bd
A symmetric curvature function $F\in C^{2,\alpha}(\C_+)\ii C^0(\bar \C_+)$
positively homogeneous of degree $d_0>0$ is said to be of class $(K)$ if
\begin{equation}\lae{1.27}
F_i=\pd F{\kappa}i>0\q \text{in } \C_+,
\end{equation}
\begin{equation}\lae{1.28}
\fmo{\pa \C_+}=0,
\end{equation}
and
\begin{equation}\lae{1.29}
F^{ij,kl}\h_{ij}\h_{kl}\le F^{-1}(F^{ij}\h_{ij})^2-F^{ik}\tilde
h^{jl}\h_{ij}\h_{kl}\qq\A\,\h\in\msc S,
\end{equation}
or, equivalently, if we set $\hat F=\log F$,
\begin{equation}\lae{1.30}
\hat F^{ij,kl}\h_{ij}\h_{kl}\le -\hat F^{ik}\tilde
h^{jl}\h_{ij}\h_{kl}\qq\A\,\h\in\msc S,
\end{equation}
where $F$ is evaluated at $(h_{ij})$.
\ed
\br\lar{1.4}\mbox{}
\bir
\item The main difference in the new definition is that we no longer
assume $F$ to be concave. Instead, we deduce from \re{1.30} that $\hat F=\log F$
is concave, which is sufficient to apply the higher regularity results once the
$C^2$-estimates are established.

\item We conclude immediately that products of functions of class $(K)$ stay in
this class, as is the case for positive powers.

\item If one wants to prove that a particular function is of class $(K)$ it might
be helpful to verify the formally less restrictive inequality \re{1.11} instead of
\re{1.29}.
\ei
\er

We immediately deduce from \re{1.29}

\bl\lal{1.2}
Let $F$ be of class $(K)$, let $\kappa_r$ be the largest eigenvalue of $(h_{ij})\in\msc
S_+$, then, for any $(\h_{ij})\in\msc S$ we have
\begin{equation}\lae{1.31}
F^{ij,kl}\h_{ij}\h_{kl}\le F^{-1}(F^{ij}\h_{ij})^2-\kappa_r^{-1}F^{ij}\h_{ir}\h_{jr},
\end{equation}
where $F$ is evaluated at $(h_{ij})$.
\el

Let $H_k$ be the symmetric polynomial of order $k$
\begin{equation}
H_k(\kappa_i)=\sum_{i_1<\dotsb<i_k}\kappa_{i_1}\dotsb\kappa_{i_k},\q 1\le k\le n,
\end{equation}
$\sigma_k=(H_k)^{1/k}$ and $\tilde \sigma_k$ the inverses of $\sigma_k$
\begin{equation}
\tilde \sigma_k(\kappa_i)=\frac{1}{\sigma_k(\kappa_i^{-1})} \,\raise 2pt \hbox{,}
\end{equation}
then, we proved in \ci[Lemma 1.5]{cg4}, see also \ci{cg6}, that the $\tilde \sigma_k$
are of class $(K)$.

Unfortunately, the class $(K)$ is too large to prove existence results in the
Lorentzian case. Instead, we have to consider a subclass $(K^*)$ which is defined
by the additional technical assumption
\bd\lad{1.3}
A function $F\in (K)$ is said to be of class $(K^*)$ if there exists
$0<\e_0=\e_0(F)$ such that

\begin{equation}\lae{1.15}
\e_0 F H\le F^{ij} h_{ik}h^k_j\,,
\end{equation}
for any $(h_{ij})\in \msc S_+$, where $F$ is evaluated at $(h_{ij})$. $H$ represents
the mean curvature, i.e. the trace of $(h_{ij})$.
\ed
\noindent Here, the index is raised with respect to the Euclidean
metric.

Evidently, $F=\sigma_n=\tilde \sigma_n$ is of class $(K^*)$ since
\begin{equation}
F^{ij}=\frac {1}{n}F\tilde h^{ij}\,,
\end{equation}
where $(\tilde h^{ij})=(h_{ij})^{-1}$.

On the other hand, the $\tilde \sigma_k,1\le k<n$, do not seem to  belong to $(K^*)$ as
is easily checked for $k=1$, while their inverses, the $\sigma_k$, fulfill \re{1.15}.
However, we shall show in \rp{1.9} below that functions of the form $FK$, where
$F\in (K)$ and $K=\sigma_n$ belong to $(K^*)$.

We should note that any symmetric $F\in C^1(\C_+)$, positively homogeneous of
degree $d_0$, with $F_i>0$ satisfies the estimate

\begin{equation}
F^{ij}h_{ik}h_j^k\le d_0FH
\end{equation}
for any $(h_{ij})\in\msc S_+$.

Before we  establish some properties of $(K^*)$, we need the following
definition.

\bd
A symmetric curvature function $F\in C^{2,\alpha}(\C_+)$
positively homogeneous of degree $d_0>0$ is said to be of class $(K_b)$, if it
satisfies the conditions of a function in class $(K)$ except the relation \re{1.28}. 
\ed

\bl\lal{1.8}
Any $F\in (K_b)$ is   bounded on bounded subsets of $\C_+$ and positive.
\el

\bp
First, we note that $F>0$ because of the homogeneity and Euler's formula. Let
$\hat F=\log F$ and consider $\kappa=(\kappa^i)\in \C_+$; in view of the concavity of
$\hat F$ we deduce
\begin{equation}
\hat F(\kappa)\le \hat F(1,\dots,1)+\hat F_i(1,\dots,1) (\kappa^i-1),
\end{equation}
i.e. $\hat F$ is locally bounded from above.
\ep

Now, we can prove
\bpp \lap{1.9} \mbox{}
\bir
\item Let $F\in (K^*)$ and $r>0$, then $F^r\in (K^*)$.
\item Let $F\in (K_b)$ and $K\in (K^*)$, then $FK \in (K^*)$.
\item The $F\in (K)$ satisfying 
\begin{equation}\lae{1.37}
F_i\kappa_i\ge \e_0F\q \forall\, i,
\end{equation}
with some positive $\e_0=\e_0(F)$, are of class $(K^*)$, and they are precisely
those, that can  be written in the form
\begin{equation}\lae{1.38}
F=GK^a,\q a>0,
\end{equation}
where $G\in (K_b)$ and  $K=\sigma_n$.
\item If $n=2$, any $F\in (K^*)$ satisfies \re{1.37}, i.e. the functions in
$(K^*)$ are exactly those given in \re{1.38}.
\ei
\epp

\bp The demonstration  of the first two properties is straight-forward, since
the product $FK$, where $F\in(K_b)$ and $K\in (K)$, can be extended as a
continuous function to $\bar\C_+$ vanishing on the boundary, so that $FK\in (K)$.

To prove (iii), we first note that any $F\in (K)$ satisfying \re{1.37} certainly
belongs to $(K^*)$, and for any $F$ of the form \re{1.38} the preceding estimate
is valid. Thus, let us assume that $F\in (K^*)$ is given for which \re{1.37} holds.
Let $\e >0$ and set $G=FK^{-\e}$. We shall show  $G\in (K_b)$, if $\e$ is small,
completing the proof of (iii).

As before, indicate the logarithm of a function  by a hat; then
\begin{equation}
\hat G_i=\hat F_i-\e\hat K_i\ge (\e_0-\tfrac{\e}{n})\kappa_i^{-1}>0,
\end{equation}
if $\e<n\e_0$, i.e. \re{1.27} is satisfied.

The inequality \re{1.30} is valid, because
this inequality becomes an equality when evaluated with  $\hat F=\hat K$.

Finally, let us derive property (iv). Assume $n=2$, and let $F\in (K^*)$, which,
without loss of generality, should be homogeneous of degree 1. Consider
$\kappa=(\kappa^1,\kappa^2)\in \C_+$ and suppose for simplicity that $\kappa^1\le \kappa^2$, then
\begin{equation}
F_2\kappa^2\le F_1\kappa^1,
\end{equation}
cf. \re{1.14}, and
\begin{equation}\lae{1.42}
F=F_1\kappa^1+F_2\kappa^2.
\end{equation}
Suppose that there is a sequence $\kappa_\e$, with $\kappa_\e^1\le\kappa_\e^2$, such that
$\hat F_2\kappa_\e^2$ tends to 0. In view of the homogeneity, we may assume that
\begin{equation}
H=\kappa_\e^1+\kappa_\e^2=1,
\end{equation}
so that we conclude from \re{1.15} and \re{1.42}
\begin{equation}
\e_0\le (\hat F_1\kappa_\e^1)\kappa_\e^1+(\hat F_2\kappa_\e^2)\kappa_\e^2\le
\kappa_\e^1+\frac{\e_0}{2},
\end{equation}
for small $\e$, i.e. $\kappa_\e^1\ge \frac{\e_0}{2}$, contradicting the assumption that
$\hat F_2(\kappa_\e)$ should tend to zero, which is only possible if
$\kappa_\e^1\rightarrow  0$.
\ep

The preceding considerations are
also applicable if the
$\kappa_i$ are the principal curvatures of a hypersurface $M$ with metric $(g_{ij})$.
$F$ can then be looked at as being defined on the space of all symmetric tensors
$(h_{ij})$ with eigenvalues $\kappa_i$ with respect to the metric.
\begin{equation}
F^{ij}=\pdc Fh{{ij}}
\end{equation}
is then a contravariant tensor of second order. Sometimes it will be convenient
to circumvent the dependence on the metric by considering $F$ to depend on the
mixed tensor
\begin{equation}
h_j^i=g^{ik}h_{kj}.
\end{equation}
Then,
\begin{equation}
F_i^j=\pdm Fhji
\end{equation}
is also a mixed tensor with contravariant index $j$ and covariant index $i$.

\section{Notations and preliminary results}\las 2

The main objective of this section is to state the equations of Gau{\ss}, Codazzi,
and Weingarten for hypersurfaces. In view of the subtle but important
difference  that is to be seen in the \tit{Gau{\ss} equation} depending on the
nature of the ambient space---Riemannian or Lorentzian---, which we already
mentioned in the introduction, we shall formulate the governing equations of a
hypersurface $M$ in a pseudo-riemannian \di{(n+1)} space $N$, which is either
Riemannian or Lorentzian. Geometric quantities in $N$ will be denoted by
$(\bar g_{\alpha\beta}),(\riema \alpha\beta\gamma\delta)$, etc., and those in $M$ by $(g_{ij}), (\riem
ijkl)$, etc. Greek indices range from $0$ to $n$ and Latin from $1$ to $n$; the
summation convention is always used. Generic coordinate systems in $N$ resp.
$M$ will be denoted by $(x^\alpha)$ resp. $(\x^i)$. Covariant differentiation will
simply be indicated by indices, only in case of possible ambiguity they will be
preceded by a semicolon, i.e. for a function $u$ in $N$, $(u_\alpha)$ will be the
gradient and
$(u_{\alpha\beta})$ the Hessian, but e.g., the covariant derivative of the curvature
tensor will be abbreviated by $\riema \alpha\beta\gamma{\delta;\e}$. We also point out that
\begin{equation}
\riema \alpha\beta\gamma{\delta;i}=\riema \alpha\beta\gamma{\delta;\e}x_i^\e
\end{equation}
with obvious generalizations to other quantities.

Let $M$ be a \tit{space-like} hypersurface, i.e. the induced metric is Riemannian,
with a differentiable normal $\n$. We define the signature of $\n$, $\sigma=\s(\n)$, by
\begin{equation}
\sigma=\bar g_{\alpha\beta}\n^\alpha\n^\beta=\spd \n\n.
\end{equation}
In case $N$ is Lorentzian, $\sigma=-1$, and $\n$ is time-like.

In local coordinates, $(x^\alpha)$ and $(\x^i)$, the geometric quantities of the
space-like hypersurface $M$ are connected through the following equations
\begin{equation}\lae{2.3}
x_{ij}^\alpha=-\sigma h_{ij}\n^\alpha
\end{equation}
the so-called \tit{Gau{\ss} formula}. Here, and also in the sequel, a covariant
derivative is always a \tit{full} tensor, i.e.

\begin{equation}
x_{ij}^\alpha=x_{,ij}^\alpha-\ch ijk x_k^\alpha+\cha \beta\gamma\alpha x_i^\beta x_j^\gamma.
\end{equation}
The comma indicates ordinary partial derivatives.

In this implicit definition the \tit{second fundamental form} $(h_{ij})$ is taken
with respect to $-\sigma\n$.

The second equation is the \tit{Weingarten equation}
\begin{equation}
\n_i^\alpha=h_i^k x_k^\alpha,
\end{equation}
where we remember that $\n_i^\alpha$ is a full tensor.

Finally, we have the \tit{Codazzi equation}
\begin{equation}
h_{ij;k}-h_{ik;j}=\riema\alpha\beta\gamma\delta\n^\alpha x_i^\beta x_j^\gamma x_k^\delta
\end{equation}
and the \tit{Gau{\ss} equation}
\begin{equation}
\riem ijkl=\sigma \{h_{ik}h_{jl}-h_{il}h_{jk}\} + \riema \alpha\beta\gamma\delta x_i^\alpha x_j^\beta x_k^\gamma
x_l^\delta.
\end{equation}
Here, the signature of $\n$ comes into play.

Now, let us assume that $N$ is a globally hyperbolic Lorentzian manifold with a
\tit{compact} Cauchy surface. As we have already pointed out in the introduction,
$N$ is then a topological product $\R[]\times \mc S_0$, where $\mc S_0$ is a
compact Riemannian manifold, and there exists a Gaussian coordinate system
$(x^\alpha)$, such that the metric in $N$ has the form \re{0.2}. We also assume that
the coordinate system is \tit{future oriented}, i.e. the time coordinate $x^0$
increases on future directed curves. Hence, the \tit{contravariant} time-like
vector $(\x^\alpha)=(1,0,\dotsc,0)$ is future directed as is its \tit{covariant}
version
$(\x_\alpha)=e^{2\psi}(-1,0,\dotsc,0)$.

Let $M=\graph \fv u\so$ be a space-like hypersurface
\begin{equation}
M=\set{(x^0,x)}{x^0=u(x),\,x\in\mc S_0},
\end{equation}
then the induced metric has the form
\begin{equation}
g_{ij}=e^{2\psi}\{-u_iu_j+\sigma_{ij}\}
\end{equation}
where $\sigma_{ij}$ is evaluated at $(u,x)$, and its inverse $(g^{ij})=(g_{ij})^{-1}$ can
be expressed as
\begin{equation}\lae{2.10}
g^{ij}=e^{-2\psi}\{\sigma^{ij}+\frac{u^i}{v}\frac{u^j}{v}\},
\end{equation}
where $(\sigma^{ij})=(\sigma_{ij})^{-1}$ and
\begin{equation}\lae{2.11}
\begin{aligned}
u^i&=\sigma^{ij}u_j\\
v^2&=1-\sigma^{ij}u_iu_j\equiv 1-\abs{Du}^2.
\end{aligned}
\end{equation}
Hence, $\graph u$ is space-like if and only if $\abs{Du}<1$.

The covariant form of a normal vector of a graph looks like
\begin{equation}
(\n_\alpha)=\pm v^{-1}e^{\psi}(1, -u_i).
\end{equation}
and the contravariant version is
\begin{equation}
(\n^\alpha)=\mp v^{-1}e^{-\psi}(1, u^i).
\end{equation}
Thus, we have
\br Let $M$ be space-like graph in a future oriented coordinate system. Then, the
contravariant future directed normal vector has the form
\begin{equation}
(\n^\alpha)=v^{-1}e^{-\psi}(1, u^i)
\end{equation}
and the past directed
\begin{equation}\lae{2.15}
(\n^\alpha)=-v^{-1}e^{-\psi}(1, u^i).
\end{equation}
\er

In the Gau{\ss} formula \re{2.3} we are free to choose the future or past directed
normal, but we stipulate that we always use the past directed normal for reasons
that will be apparent in a moment.

Look at the component $\alpha=0$ in \re{2.3} and obtain in view of \re{2.15}

\begin{equation}\lae{2.16}
e^{-\psi}v^{-1}h_{ij}=-u_{ij}-\cha 000\mspace{1mu}u_iu_j-\cha 0i0
\mspace{1mu}u_j-\cha 0j0\mspace{1mu}u_i-\cha ij0.
\end{equation}
Here, the covariant derivatives a taken with respect to the induced metric of
$M$, and
\begin{equation}
-\cha ij0=e^{-\psi}\bar h_{ij},
\end{equation}
where $(\bar h_{ij})$ is the second fundamental form of the hypersurfaces
$\{x^0=\const\}$.

An easy calculation shows
\begin{equation}
\bar h_{ij}e^{-\psi}=-\tfrac{1}{2}\dot\sigma_{ij} -\dot\psi\sigma_{ij},
\end{equation}
where the dot indicates differentiation with respect to $x^0$.

Let us assume for the moment that the Gaussian coordinate system is
\tit{normal}, i.e. $\psi\equiv 0$, then
\begin{equation}
\bar h_{ij}=-\tfrac{1}{2}\dot\sigma_{ij},
\end{equation}
and the mean curvature of the level hypersurfaces, $\bar H=\sigma^{ij}\bar h_{ij}$,
satisfies the equation
\begin{equation}
\Dot{\bar H}=\bar R_{\alpha\beta}\n^\alpha\n^\beta+{\bar h_{ij}\bar h^{ij},}
\end{equation}
as one can easily check. If we assume now, that the \tit{time-like convergence}
condition holds in $N$, i.e
\begin{equation}
\bar R_{\alpha\beta}\x^\alpha\x^\beta\ge 0
\end{equation}
for all time-like $(\x^\alpha)$, then we deduce that $\bar H$ is monotone increasing
in time.

Thus, we see that our intuitive understanding, namely, that lower barriers, as
defined in \rd{0.1}, should lie in the past of upper barriers is generically in
accordance with Lorentzian geometry if we evaluate the second fundamental
form with respect to the past directed normal.

\bd
A closed, space-like hypersurface $M$ is said to be \tit{convex} (\tit{strictly
convex}) if its second fundamental form evaluated with respect to the past
directed normal is \tit{positive} semi-definite (definite).
\ed

\br
If in a particular setting the second fundamental forms of the barriers involved
are negative semi-definite, when evaluated with respect to the past directed
normal, then, changing the roles of the future and past directed light cones will
establish the preferred situation, where convexity means non-negative principal
curvatures.
\er

Next, let us analyze under which condition a space-like hypersurface $M$ can be
written as a graph over the Cauchy hypersurface $\mc S_0$.

We first need
\bd

Let $M$ be a  closed, space-like hypersurface in $N$. Then,
\bi[(i)]
\item
$M$ is said to be \tit{achronal}, if no two points in $M$ can be connected by a
future directed time-like curve.

\item
$M$ is said to \tit{separate} $N$, if $N\sminus M$ is disconnected.
\ei
\ed

We can now prove

\bpp\lap{2.5}
Let $N$ be connected and globally hyperbolic, $\mc S_0
\su N$ a compact Cauchy hypersurface, and $M\su N$ a compact,
connected space-like hypersurface of class $C^m, m\ge 1$. Then, $M=\graph \fv
u\so$ with
$u\in C^m({\mc S}_0)$ iff $M$ is achronal.
\epp
\bp
(i) We first show that an achronal $M$ is a graph over $\so$. Let $(x^\alpha)$ be 
the special coordinate system associated with $\mc S_0$ such that
$\so=\set {p\in N}{x^0(p)=0}$, and let $p \in M$ be arbitrary, $p=\left(x^0(p),
x(p)\right)$. Since $M$ is achronal, the time-like curve $\{\gamma_p\}=\set
{\left(x^0,x(p)\right)}{x^0\in
\R[]}$ through $\left(0,x(p)\right)\in \so$ intersects $M$ exactly once, and we conclude
that $M=\graph \fv uG$ with $u \in C^0(G)$, where $G\su\so$ is closed. But $G$ is
also open, and hence $G=\so$, for otherwise, there would be $q\in M$ such that
$\dot \gamma_q\in T_q(M)$, which is impossible since $M$ has a continuous time-like
normal.

Furthermore, there exists a neighbourhood $\mc U=\mc U(p)$ in $N$
and a function $\Phi\in C^m(\mc U)$ with time-like gradient such that
\begin{equation}
\mc U\ii M=\set {(x^0,x)}{\Phi(x^0,x)=0}.
\end{equation}
$M$ is connected with a continuous time-like normal. Thus, we obtain
\begin{equation}
\pd {\Phi}x0=\spd {D\Phi}{\pd {}x0}\ne0,
\end{equation}
and we deduce from the implicit function theorem, that there is a neighbourhood
$\mc V$ of $x(p)$ in $\so$ and a possibly smaller neighbourhood $\widetilde {\mc
U}$ of $p$ such that
\begin{equation}
\widetilde {\mc U}\ii M=\graph \fv\f{\mc V}\,,\q \f\in C^m(\mc V).
\end{equation}
Hence, $\f=\fv u{\mc V}$ and $u$ is of class $C^m$.

(ii) To demonstrate the reverse implication, we use the fact that $M$ is achronal
if
$M$ separates $N$, cf. \ci[p. 427]{bon}, and observe that any graph over $\so$
separates $N$.
\ep

In \ci[p. 427]{bon} it is also proved that a closed, connected, space-like
hypersurface
$M$ is achronal if $N$ is simply connected. Hence, we infer

\br
Assume that the Cauchy hypersurface $\so$ is homeomorphic to $S^n, n\ge 2$,
then any closed, connected space-like hypersurface $M$ is a graph over $\so$.
\er

One of the assumptions in \rt{0.2} is that there exists a strictly convex function
$\chi\in C^2(\bar \Omega)$. We shall state sufficient geometric conditions
guaranteeing the existence of such a function.

\bl
Let $N$ be globally hyperbolic, $\so$ a Cauchy hypersurface, $(x^\alpha)$ a special
coordinate system associated with $\so$, and $\bar \Omega\su N$ be compact. Then,
there exists a strictly convex function $\chi\in C^2(\bar \Omega)$ provided the level
hypersurfaces $\{x^0=\const\}$ that intersect $\bar \Omega$ are strictly convex.
\el

\bp
For greater clarity set $t=x^0$, i.e. $t$ is a globally defined time function. Let
$x=x(\x)$ be a local representation for $\{t=\const\}$, and $t_i,t_{ij}$ be the
covariant derivatives of $t$ with respect to the induced metric, and
$t_\alpha,t_{\alpha\beta}$ be the covariant derivatives in $N$, then
\begin{equation}
0=t_{ij}=t_{\alpha\beta}x_i^\alpha x_j^\beta+t_\alpha x_{ij}^\alpha,
\end{equation}
and therefore,
\begin{equation}\lae{2.26}
t_{\alpha\beta}x_i^\alpha x_j^\beta=-t_\alpha x_{ij}^\alpha=-\bar h_{ij} t_\alpha \n^\alpha.
\end{equation}
Here, $(\n^\alpha)$ is past directed, i.e. the right-hand side in \re{2.26} is positive
definite in $\bar \Omega$, since $(t_\alpha)$ is also past directed.

Choose $\lambda>0$ and define $\chi=e^{\lambda t}$, so that
\begin{equation}
\chi_{\alpha\beta}=\lambda^2e^{\lambda t} t_\alpha t_\beta+\lambda e^{\lambda t} t_{\alpha\beta}.
\end{equation}

Let $p\in \Omega$ be arbitrary, $\mc S=\{t=t(p)\}$ be the level hypersurface through
$p$, and $(\h^\alpha)\in T_p(N)$. Then, we conclude
\begin{equation}
e^{-\lambda t}\chi_{\alpha\beta}\h^\alpha \h^\beta=\lambda^2\abs{\h^0}^2+\lambda t_{ij}\h^i\h^j+2\lambda t_{0j}\h^0\h^i,
\end{equation}
where $t_{ij}$ now represents the left-hand side in  \re{2.26}, and we infer further
\begin{equation}
\begin{aligned}
e^{-\lambda t}\chi_{\alpha\beta}\h^\alpha\h^\beta&\ge \tfrac{1}{2} \lambda^2 {\abs\h^0}^2 +[\lambda
\e-c_\e]\sigma_{ij}\h^i\h^j\\
&\ge \tfrac{\e}{2}\lambda \{-\abs{\h^0}^2+\sigma_{ij}\h^i\h^j\}
\end{aligned}
\end{equation}
for some $\e>0$, and where $\lambda$ is supposed to be large. Therefore, we
have in $\bar \Omega$
\begin{equation}
\chi_{\alpha\beta}\ge c\bar g_{\alpha\beta}\,,\q c>0,
\end{equation}
i.e. $\chi$ is strictly convex.
\ep

Sometimes, we need a Riemannian reference metric, e.g. if we want to estimate
tensors. Since the Lorentzian metric can be expressed as
\begin{equation}
\bar g_{\alpha\beta}dx^\alpha dx^\beta=e^{2\psi}\{-{dx^0}^2+\sigma_{ij}dx^i dx^j\},
\end{equation}
we define a Riemannian reference metric $(\tilde g_{\alpha\beta})$ by
\begin{equation}
\tilde g_{\alpha\beta}dx^\alpha dx^\beta=e^{2\psi}\{{dx^0}^2+\sigma_{ij}dx^i dx^j\}
\end{equation}
and we abbreviate the corresponding norm of a vectorfield $\h$ by
\begin{equation}
\nnorm \h=(\tilde g_{\alpha\beta}\h^\alpha\h^\beta)^{1/2},
\end{equation}
with similar notations for higher order tensors.

\section{The evolution problem}\las 3

Solving the problem \re{0.1} consists of two steps: first, one has to prove a
priori estimates, and secondly, one has to find a procedure which, with the help
of the priori estimates, leads to a solution of the problem.

When we first considered the problem for $F\in (K)$ in the Riemannian case, we
used an evolutionary approach, which was rather aesthetic but had the
short-coming that for technical reasons the sectional curvatures of the ambient
space had to be non-positive, cf. \ci{cg2}. We were able to overcome this
technical obstruction in \ci{cg4}, where we used the method of successive
approximation to prove  existence. An important ingredient of that proof was
the property of the class $(K)$ to be closed under \tit{elliptic regularization},
see \ci[Section 1]{cg4} for details. However, the subclass $(K^*)$ is not closed
under elliptic regularization, so that this method of proof fails in the Lorentzian
case. But, fortunately, we can apply the evolutionary approach without making
any sacrifices with respect to the sectional curvatures of the ambient space,
since the unfavourable sign condition that forces us to consider the class $(K^*)$
instead of $(K)$ eliminates that particular technical obstruction.

For greater transparency, we look at the problem in a pseudo-rieman\-nian space
$N$, where, as already stated in \rs{2}, we, really, only have the Riemannian and
the Lorentzian case in mind. Properties like space-like, achronal, etc., however, 
only make sense, when $N$ is Lorentzian and should be ignored otherwise.

We want to prove that the equation 
\begin{equation}
F=f
\end{equation}
has a solution. For technical reasons, it is convenient to solve instead the
equivalent equation
\begin{equation}\lae{3.2}
\F(F)=\F(f),
\end{equation}
where $\F$ is a real function defined on $\R[]_+$ such that
\begin{equation}
\dot\F>0\q \tup{and}\q \ddot\F\le 0.
\end{equation}

For notational reasons, let us abbreviate
\begin{equation}
\tilde f=\F(f).
\end{equation}

We also point out that we may---and shall---assume without loss of generality
that $F$ is homogeneous of degree 1.

To solve \re{3.2} we look at the evolution problem
\begin{equation}\lae{3.5}
\begin{aligned}
\dot x&=-\s(\F-\tilde f)\n,\\
x(0)&=x_0,
\end{aligned}
\end{equation}
where $x_0$ is an embedding of an initial strictly convex, compact, space-like
hypersurface $M_0$, $\F=\F(F)$, and $F$ is evaluated at the principal curvatures
of the flow hypersurfaces $M(t)$, or, equivalently, we may assume that $F$
depends on the second fundamental form $(h_{ij})$ and the metric $(g_{ij})$ of
$M(t)$; $x(t)$ is the embedding of $M(t)$ and $\sigma$ the signature of the (past
directed) normal $\n=\n(t)$.

This is a parabolic problem, so short-time existence is guaranteed---the proof
in the Lorentzian case is identical to that in the Riemannian case, cf. \ci[p.
622]{cg2}---, and under suitable assumptions, we shall be able to prove that the
solution exists for all time and converges to a stationary solution if $t$ goes to
infinity.

There is a slight ambiguity in the notation, since we also call the
evolution parameter \tit{time}, but this lapse shouldn't cause any
misunderstandings.

Next, we want to show how the metric, the second fundamental form, and the
normal vector of the hypersurfaces $M(t)$ evolve. All time derivatives are
\tit{total} derivatives. The proofs are identical to those of the corresponding
results in a Riemannian setting, cf. \ci[Section 3]{cg2}, and will be omitted.

\bl[Evolution of the metric]
The metric $g_{ij}$ of $M(t)$ satisfies the evolution equation
\begin{equation}
\dot g_{ij}=-2\s(\F-\tilde f)h_{ij}.
\end{equation}
\el

\bl[Evolution of the normal]
The normal vector evolves according to
\begin{equation}\lae{3.7}
\dot \n=\nabla_M(\F-\tilde f)=g^{ij}(\F-\tilde f)_i x_j.
\end{equation}
\el

\bl[Evolution of the second fundamental form]
The second fundamental form evolves according to
\begin{equation}\lae{3.8}
\dot h_i^j=(\F-\tilde f)_i^j+\sigma (\F-\tilde f) h_i^k h_k^j+\sigma (\F-\tilde f) \riema
\alpha\beta\gamma\delta\n^\alpha x_i^\beta \n^\gamma x_k^\delta g^{kj}
\end{equation}
and
\begin{equation}
\dot h_{ij}=(\F-\tilde f)_{ij}-\sigma (\F-\tilde f) h_i^k h_{kj}+\sigma (\F-\tilde f) \riema
\alpha\beta\gamma\delta\n^\alpha x_i^\beta \n^\gamma x_j^\delta.
\end{equation}
\el

\bl[Evolution of $(\F-\tilde f)$]
The term $(\F-\tilde f)$ evolves according to the equation
\begin{equation}\lae{3.10}
\begin{aligned}
{(\F-\tilde f)}^\prime-\dot\F F^{ij}(\F-\tilde f)_{ij}=&\msp[3]\sigma \dot \F
F^{ij}h_{ik}h_j^k (\F-\tilde f)\\
 &+\sigma\tilde f_\alpha\n^\alpha (\F-\tilde f)\\
&+\sigma\dot\F F^{ij}\riema \alpha\beta\gamma\delta\n^\alpha x_i^\beta \n^\gamma x_j^\delta (\F-\tilde f),
\end{aligned}
\end{equation}
where
\begin{equation}
(\F-\tilde f)^{\prime}=\frac{d}{dt}(\F-\tilde f)
\end{equation}
and
\begin{equation}
\dot\F=\frac{d}{dr}\F(r).
\end{equation}
\el

From \re{3.8} we deduce with the help of the Ricci identities a parabolic equation
for the second fundamental form
\bl
The mixed tensor $h_i^j$ satisfies the parabolic equation
\begin{equation}\raisetag{-78pt}\lae{3.13}
\begin{aligned}
\dot h_i^j-\dot\F F^{kl}h_{i;kl}^j&=\sigma \dot\F F^{kl}h_{rk}h_l^rh_i^j-\sigma\dot\F F
h_{ri}h^{rj}+\sigma (\F-\tilde f) h_i^kh_k^j\\
&\hp{+}-\tilde f_{\alpha\beta} x_i^\alpha x_k^\beta g^{kj}+\sigma \tilde f_\alpha\n^\alpha h_i^j+\dot\F
F^{kl,rs}h_{kl;i}h_{rs;}^{\hphantom{rs;}j}\\
&\hp{=}+\ddot \F F_i F^j+2\dot \F F^{kl}\riema \alpha\beta\gamma\delta x_m^\alpha x_i ^\beta x_k^\gamma
x_r^\delta h_l^m g^{rj}\\
&\hp{=}-\dot\F F^{kl}\riema \alpha\beta\gamma\delta x_m^\alpha x_k ^\beta x_r^\gamma x_l^\delta
h_i^m g^{rj}-\dot\F F^{kl}\riema \alpha\beta\gamma\delta x_m^\alpha x_k ^\beta x_i^\gamma x_l^\delta h^{mj} \\
&\hp{=}+\sigma\dot\F F^{kl}\riema \alpha\beta\gamma\delta\n^\alpha x_k^\beta\n^\gamma x_l^\delta h_i^j-\sigma\dot\F F
\riema \alpha\beta\gamma\delta\n^\alpha x_i^\beta\n^\gamma x_m^\delta g^{mj}\\
&\hp{=}+\sigma (\F-\tilde f)\riema \alpha\beta\gamma\delta\n^\alpha x_i^\beta\n^\gamma x_m^\delta g^{mj}\\
&\hp{=}+\dot\F F^{kl}\bar R_{\alpha\beta\gamma\delta;\e}\{\n^\alpha x_k^\beta x_l^\gamma x_i^\delta
x_m^\e g^{mj}+\n^\alpha x_i^\beta x_k^\gamma x_m^\delta x_l^\e g^{mj}\}.
\end{aligned}
\end{equation}
\el

The proof is identical to that of the corresponding result in the Riemannian case,
cf. \ci[Lemma 7.1 and Lemma 7.2]{cg2}; we only have to keep track of the
signature of the normal in the more general pseudo-riemannian setting. 

If we had assumed $F$ to be homogeneous of degree $d_0$ instead of 1, then, we
would have to replace the explicit term $F$---occurring twice in the preceding
lemma---by $d_0F$.

We also point out that the technical differences we encounter, due to the nature
of the ambient space---Riemannian or Lorentzian---, stem from the
alternating sign of $\sigma$ in \re{3.13}.

\br\lar{3.6}
In view of the maximum principle, we immediately deduce from \re{3.10} that the
term $(\F-\tilde f)$ has a sign during the evolution if it has one at the beginning,
e.g., if the starting hypersurface $M_0$ is the upper barrier $M_2$, then
$(\F-\tilde f)$ is non-negative, or equivalently,
\begin{equation}\lae{3.14}
F\ge f.
\end{equation}
\er

\section{Lower order estimates}\las 4

From now on, we stick to our original assumption that the ambient space is
globally hyperbolic with a compact Cauchy hypersurface $\so$. The barriers
$M_i$ are then graphs over $\so, M_i=\graph u_i$, because they are achronal, cf.
\rp{2.5}, and we have
\begin{equation}\lae{4.1}
u_1\le u_2,
\end{equation}
for $M_1$ should lie in the past of $M_2$, and the enclosed domain is supposed to
be connected. Moreover, in view of the Harnack inequality, the strict inequality
is valid in \re{4.1} unless the barriers coincide and are a solution to our
problem, cf.  the proof of \rl{4.1}.

Let us look at the evolution equation \re{3.5} with initial hypersurface $M_0$
equal to $M_2$. Then, because of the short-time existence, the
evolution will exist on a maximal time interval
$I=[0,T^*),T^*\le
\un$, as long as the evolving hypersurfaces are space-like, strictly convex and
smooth.

Furthermore, since the initial hypersurface is a graph over $\so$, we can write
\begin{equation}
M(t)=\graph\fu{u(t)}S0\q \A\,t\in I,
\end{equation}
where $u$ is defined in the cylinder $Q_{T^*}=I\times \so$. We then deduce from
\re{3.5}, looking at the component $\alpha=0$, that $u$ satisfies a parabolic equation
of the form
\begin{equation}\lae{4.3}
\dot u=-e^{-\psi}v^{-1}(\F-\tilde f),
\end{equation}
where we  use the notations in \rs{2}, and where we emphasize that the time
derivative is a total derivative, i.e.
\begin{equation}\lae{4.4}
\dot u=\pde ut+u_i\dot x^i.
\end{equation}

Since the past directed normal can be expressed as
\begin{equation}
(\n^\alpha)=-e^{-\psi}v^{-1}(1,u^i),
\end{equation}
we conclude from \re{3.5}, \re{4.3}, and \re{4.4}
\begin{equation}
\pde ut=-e^{-\psi}v(\F-\tilde f).
\end{equation}
Thus, $\pde ut$ is non-positive in view of \rr{3.6}.

Next, let us state our first a priori estimate

\bl\lal{4.1}
During the evolution the flow hypersurfaces stay in $\bar \Omega$.
\el

\bp
Since  $\pde ut$ is non-positive, we only have to consider the case that the
flow reaches the boundary component $M_1$.   Suppose that  the flow
hypersurfaces would touch
$M_1$ for the first time at time $t=t_0$ in $x_0\in M_1$, then, we deduce from
the equation \re{2.16} and the maximum principle, that $x_0\in \varSigma$ and
conclude further that, in view of the relation
\re{3.14}, the  Harnack inequality can be applied to $(u-u_1)$
to yield 
$M(t_0)=M_1$, and hence, that $M_1$  is already a
solution to our problem; the flow would become stationary for $t\ge t_0$.
\ep

\br\lar{4.2}
It is important to allow non-convex lower barriers, because the \tit{big bang}
and
\tit{big crunch} hypotheses of the standard cosmological model assert that
there are sequences $M_{1,k}$ and
$M_{2,k}$ of closed, achronal, space-like hypersurfaces such that, in our setting,
$M_{i,k}=\graph \fv {u_{i,k}}{\so}$, for $i=1,2$, 
\begin{equation}
\lim_{k\rightarrow \un} \sup_\so u_{1,k}=-\un,\q \lim_{k\rightarrow \un}
\inf_\so u_{2,k}=\un,
\end{equation}
and the principal curvatures with respect to the past directed normal of
$M_{1,k}$ tend to $-\un$, while those of $M_{2,k}$ tend to $\un$.

Thus, the $M_{2,k}$ could serve as upper barriers for our purposes, but the
$M_{1,k}$ would fail to be lower barriers, if we would only consider 
convex hypersurfaces. 
\er

As a consequence of \rl{4.1} we obtain
\begin{equation}
\inf_{\so} u_1\le u\le \sup_\so u_2\q \A\,t\in I.
\end{equation}

We are now able to derive the $C^1$-estimates, i.e. we shall show that the
hypersurfaces remain uniformly space-like, or equivalently, that the term
\begin{equation}
\tilde v=v^{-1}=\frac{1}{\sql}
\end{equation}
is uniformly bounded.

In the Riemannian case, $C^1$-estimates for closed, convex hypersurfaces can
only be derived if they are graphs in a \tit{normal} Gaussian coordinate system,
in the Lorentzian case the Gaussian coordinate system no longer needs to be
normal, and also, the convexity assumption can be relaxed to a unilateral bound
for the second fundamental form.

\bl
Let $M=\graph \fv u\so$ be a compact, space-like hypersurface represented in a
Gaussian coordinate system with  unilateral bounded  principal curvatures,
e.g.
\begin{equation}
\kappa_i\ge \kappa_0\q\A\,i.
\end{equation}
Then, the quantity $\tilde v=\frac{1}{\sql}$ can be estimated by

\begin{equation}
\tilde v\le c(\abs u,\so,\sigma_{ij},\psi,\kappa_0),
\end{equation}
where we used the notation in \re{0.2}, i.e. in the Gaussian coordinate system the
ambient metric has the form
\begin{equation}
d\bar s_N^2=e^{2\psi}\{-{dx^0}^2+\sigma_{ij}(x^0,x)dx^idx^j\}.
\end{equation}
\el
\bp
We suppose that the Gaussian coordinate system is future oriented, and that the
second fundamental form is evaluated with respect to the past-directed normal.
From formulas \re{2.10} and \re{2.11} we get
\begin{equation}
\norm{Du}^2=g^{ij}u_iu_j=e^{-2\psi}\frac{\abs{Du}^2}{v^2}\raise 2pt \hbox{,}
\end{equation}
hence, it is equivalent to find an a priori estimate for $\norm{Du}$.

Let $\lambda$ be a real parameter to be specified later, and set
\begin{equation}
w=\tfrac{1}{2}\log\norm{Du}^2+\lambda u.
\end{equation}
We may regard $w$ as being defined on $\so$; thus, there is $x_0\in\so$ such that
\begin{equation}
w(x_0)=\sup_\so w,
\end{equation}
and we conclude
\begin{equation}
0=w_i=\frac{1}{\norm{Du}^2}\,u_{ij}u^j+\lambda u_i
\end{equation}
in $x_0$, where the covariant derivatives are taken with respect to the induced
metric
$g_{ij}$, and the indices are also raised with respect to that metric.

In view of \re{2.16} we deduce further
\begin{equation}\lae{4.16}
\begin{aligned}
\lambda\norm{Du}^4&=-u_{ij}u^iu^j\\
&= e^{-\psi}\tilde vh_{ij}u^iu^j+\cha 000\msp \norm{Du}^4\\
&\hp{=}\msp[2]+2\cha 0j0\msp u^j\norm{Du}^2+\cha ij0\msp u^iu^j.
\end{aligned}
\end{equation}
Now, there holds
\begin{equation}
u^i=g^{ij}u_j=e^{-2\psi}\sigma^{ij}u_jv^{-2},
\end{equation}
and by assumption,
\begin{equation}
h_{ij}u^iu^j\ge \kappa_0\msp\norm{Du}^2,
\end{equation}
i.e. the critical terms on the right-hand side of \re{4.16} are of fourth order in
$\norm{Du}$ with bounded coefficients, and we conclude that $\norm{Du}$ can't be
too large in $x_0$ if we choose $\lambda$ such that
\begin{equation}
\lambda\le -c\msp\nnorm{\cha \alpha\beta 0}-1
\end{equation}
with a suitable constant $c$; $w$, or equivalently, $\norm{Du}$ is therefore
uniformly bounded from above.
\ep

For convex graphs over $\so$ the term $\tilde v$ is uniformly bounded as long as
they stay in a compact set. Moreover, we shall see, that $\tilde v$ satisfies a
useful parabolic equation that we shall exploit to estimate the principal
curvatures of the hypersurfaces $M(t)$ from above.

\bl[Evolution of $\tilde v$]
Consider the flow \re{3.5} in the distinguished coordinate system associated
with $\so$. Then, $\tilde v$ satisfies the evolution equation
\begin{equation}\lae{4.20}
\begin{aligned}
\dot{\tilde v}-\dot\F F^{ij}\tilde v_{ij}=&-\dot\F F^{ij}h_{ik}h_j^k\tilde v
+[(\F-\tilde f)-\dot\F F]\h_{\alpha\beta}\n^\alpha\n^\beta\\
&-2\dot\F F^{ij}h_j^k x_i^\alpha x_k^\beta \h_{\alpha\beta}-\dot\F F^{ij}\h_{\alpha\beta\gamma}x_i^\beta
x_j^\gamma\n^\alpha\\
&-\dot\F F^{ij}\riema \alpha\beta\gamma\delta\n^\alpha x_i^\beta x_k^\gamma x_j^\delta\h_\e x_l^\e g^{kl}\\
&-\tilde f_\beta x_i^\beta x_k^\alpha \h_\alpha g^{ik},
\end{aligned}
\end{equation}
where $\h$ is the covariant vector field $(\h_\alpha)=e^{\psi}(-1,0,\dotsc,0)$.
\el

\bp
We have $\tilde v=\spd \h\n$. Let $(\x^i)$ be local coordinates for $M(t)$.
Differentiating $\tilde v$ covariantly we deduce
\begin{equation}\lae{4.21}
\tilde v_i=\h_{\alpha\beta}x_i^\beta\n^\alpha+\h_\alpha\n_i^\alpha,
\end{equation}
\begin{equation}\lae{4.22}
\begin{aligned}
\tilde v_{ij}= &\msp[5]\h_{\alpha\beta\gamma}x_i^\beta x_j^\gamma\n^\alpha+\h_{\alpha\beta}x_{ij}^\beta\n^\alpha\\
&+\h_{\alpha\beta}x_i^\beta\n_j^\alpha+\h_{\alpha\beta}x_j^\beta\n_i^\alpha+\h_\alpha\n_{ij}^\alpha
\end{aligned}
\end{equation}

The time derivative of $\tilde v$ can be expressed as
\begin{equation}\lae{4.23}
\begin{aligned}
\dot{\tilde v}&=\h_{\alpha\beta}\msp\dot x^\beta\n^\alpha+\h_\alpha\dot\n^\alpha\\
&=\h_{\alpha\beta}\n^\alpha\n^\beta(\F-\tilde f)+(\F-\tilde f)^k x_k^\alpha\h_\alpha\\
&=\h_{\alpha\beta}\n^\alpha\n^\beta(\F-\tilde f)+\dot\F F^k x_k^\alpha\h_\alpha-{\tilde f}_\beta x_i^\beta
x_k^\alpha g^{ik}\h_\alpha,
\end{aligned}
\end{equation}
where we have used \re{3.7}.

Substituting \re{4.22} and \re{4.23} in \re{4.20}, and simplifying the resulting
equation with the help of the Weingarten and Codazzi equations, we arrive at the
desired conclusion.
\ep

\section{A priori estimates in the $C^2$-{norm}}\las 5

Let $M(t)$ be a solution of the evolution problem \re{3.5} with initial
hypersurface $M_0=M_2$, defined on a maximal time interval $I=[0,T^*)$. We
assume that $F$ is of class $(K^*)$ according to \rd{1.3}, homogeneous of
degree 1, and we choose
$\F(r)=\log r$; alternatively, we could  use
$\F(r)=-\frac{1}{m}r^{-m}$, $\,m\ge 1$, but with the logarithm the proof of the
$C^2$-estimates is a bit simpler. Furthermore, we suppose that there exists
a strictly convex function
$\chi\in C^2(\bar \Omega)$, i.e. there holds
\begin{equation}
\chi_{\alpha\beta}\ge c_0\bar g_{\alpha\beta}
\end{equation}
with a positive constant $c_0$.

We observe that
\begin{equation}\lae{5.2}
\begin{aligned}
\dot\chi-\dot\F F^{ij}\chi_{ij}&=[(\F-\tilde f)-\dot\F F]\chi_\alpha\n^\alpha-\dot\F
F^{ij}\chi_{\alpha\beta} x_i^\alpha x_j^\beta\\
&\le [(\F-\tilde f)-\dot\F F] \chi_\alpha\n^\alpha-c_0\dot\F F^{ij} g_{ij},
\end{aligned}
\end{equation}
where we used the homogeneity of $F$.

From \rr{3.6} we infer 
\begin{equation}\lae{5.4}
\F\ge \tilde f\qq \tup{or}\qq F\ge f,
\end{equation}
and from the results in \rs{4} that the flow stays in the compact set $\bar \Omega$.

Furthermore, due to \re{5.4} and the fact that $M_0$ is strictly convex, the $M(t)$
remain strictly convex during the evolution; hence, $\tilde v$ is uniformly
bounded.

We are now able to prove

\bl\lal{5.1}
Let $F$ be of class $(K^*)$. Then, the principal curvatures of the evolution
hypersurfaces $M(t)$ are uniformly bounded.
\el

\bp
Let $\f$ and $w$ be defined respectively by
\begin{align}
\f&=\sup\set{{h_{ij}\h^i\h^j}}{{\norm\h=1}},\\
w&=\log\f+\lambda\tilde v+\m\chi,\lae{5.6}
\end{align}
where $\lambda,\m$ are large positive parameters to be specified later. We claim that
$w$ is bounded for a suitable choice of $\lambda,\m$.

Let $0<T<T^*$, and $x_0=x_0(t_0)$, with $ 0<t_0\le T$, be a point in $M(t_0)$ such
that

\begin{equation}
\sup_{M_0}w<\sup\set {\sup_{M(t)} w}{0<t\le T}=w(x_0).
\end{equation}

We then introduce a Riemannian normal coordinate system $(\x^i)$ at $x_0\in
M(t_0)$ such that at $x_0=x(t_0,\x_0)$ we have
\begin{equation}
g_{ij}=\delta_{ij}\q \tup{and}\q \f=h_n^n.
\end{equation}

Let $\tilde \h=(\tilde \h^i)$ be the contravariant vector field defined by
\begin{equation}
\tilde \h=(0,\dotsc,0,1),
\end{equation}
and set
\begin{equation}
\tilde \f=\frac{h_{ij}\tilde \h^i\tilde \h^j}{g_{ij}\tilde \h^i\tilde \h^j}\raise 2pt
\hbox{.}
\end{equation}

$\tilde \f$ is well defined in neighbourhood of $(t_0,\x_0)$.

Now, define $\tilde w$ by replacing $\f$ by $\tilde \f$ in \re{5.6}; then, $\tilde w$
assumes its maximum at $(t_0,\x_0)$. Moreover, at $(t_0,\x_0)$ we have
\begin{equation}
\dot{\tilde \f}=\dot h_n^n,
\end{equation}
and the spatial derivatives do also coincide; in short, at $(t_0,\x_0)$ $\tilde \f$
satisfies the same differential equation \re{3.13} as $h_n^n$. For the sake of
greater clarity, let us therefore treat $h_n^n$ like a scalar and pretend that $w$
is defined by
\begin{equation}
w=\log h_n^n+\lambda\tilde v+\m\chi.
\end{equation}

At $(t_0,\x_0)$ we have $\dot w\ge 0$, and, in view of the maximum principle, we
deduce from \re{1.15}, \re{3.13}, \re{4.20}, and \re{5.2}
\begin{equation}\lae{5.13}
\begin{aligned}
0\le&\msp[3] \dot\F F h_n^n-(\F-\tilde f) h_n^n+\lambda c_1 -\lambda\e_0\dot\F F H \tilde
v\\ & +\lambda c_1[(\F-\tilde f)+\dot\F F]+\lambda c_1\dot\F F^{ij} g_{ij}\\
&+ \m c_1 [(\F-\tilde f)+\dot\F F]-\m c_0 \dot\F F^{ij} g_{ij}\\
&+\dot\F F^{ij}(\log h_n^n)_i(\log h_n^n)_j\\
&+\{\ddot\F F_n F^n +\dot\F
F^{kl,rs}h_{kl;n}h_{rs;}^{\hphantom{rs;}n}\}(h_n^n)^{-1},
\end{aligned}
\end{equation}
where we have estimated bounded terms by a constant $c_1$, assumed that
$h_n^n, \lambda$, and $\m$ are larger than $1$, and used \re{5.4} as well as the simple
observation
\begin{equation}\lae{5.14}
\abs{F^{ij}h_j^k\h_k}\le \norm\h F
\end{equation}
for any vector field $(\h_k)$, cf. \ci[Lemma 7.4]{cg2}.

Now, the last term in \re{5.13} is estimated from above by
\begin{equation}\lae{5.15}
\{\ddot\F F_n F^n+\dot\F F^{-1} F_n F^n\}(h_n^n)^{-1}-\dot \F F^{ij}
h_{in;n}h_{jn;}^{\hphantom{jn;}n}(h_n^n)^{-2},
\end{equation}
cf. \re{1.31}, where  the sum in the braces vanishes,
due to the choice of
$\F$. Moreover, because of the Codazzi equation, we have
\begin{equation}
h_{in;n}=h_{nn;i}+\riema \alpha\beta\gamma\delta\n^\alpha x_n^\beta x_i^\gamma x_n^\delta,
\end{equation}
and hence, using the  abbreviation $\bar R_i$ for the curvature term, we conclude
that \re{5.15} is bounded from above by
\begin{equation}
-(h_n^n)^{-2}\dot\F F^{ij}(h_{n;i}^n+\bar R_i)(h_{n;j}^n+\bar R_j).
\end{equation}

Thus, the terms in \re{5.13} containing the derivatives of $h_n^n$ are estimated
from above by
\begin{equation}
-2\dot\F F^{ij}(\log h_n^n)_i\bar R_j(h_n^n)^{-1}.
\end{equation}

Moreover,  $Dw$ vanishes at $\x_0$, i.e.
\begin{equation}
D\log h_n^n=-\lambda D\tilde v-\m D\chi,
\end{equation}
where only $D\tilde v$ deserves further consideration.

Replacing then $D\tilde v$ by the right-hand side of \re{4.21}, and using the
Weingarten equation and \re{5.14}, we finally conclude from \re{5.13}

\begin{equation}
\begin{aligned}
0\le &\msp[3]\dot\F F h_n^n-(\F-\tilde f)h_n^n+\lambda c_1-\lambda\e_0\dot\F
F H
\tilde v\\
&+(\lambda+\m)c_1[(\F-\tilde f)+\dot\F F]+\lambda c_1\dot\F F^{ij} g_{ij}\\
&-\m [c_0-c_1 (h_n^n)^{-1}]\dot\F F^{ij} g_{ij}
\end{aligned}
\end{equation}

Then, if we suppose $h_n^n$ to be so large that
\begin{equation}
c_1\le \tfrac{1}{2}c_0 h_n^n,
\end{equation}
and if we choose $\lambda, \m$ such that
\begin{align}\lae{5.22}
2&\le \lambda\e_0\\
\intertext{and}
4\lambda c_1&\le \m c_0
\end{align}
we derive
\begin{equation}
\begin{aligned}
0\le&\msp[2] -\tfrac{1}{2}\lambda\e_0 \dot\F F H \tilde v-(\F-\tilde f) h_n^n\\
&+(\lambda+\m)c_1 [(\F-\tilde f)+\dot\F F]+\lambda c_1.
\end{aligned}
\end{equation}

We now observe that $\dot \F F=1$, and deduce in view of \re{5.4} that $h_n^n$ is
a priori bounded at $(t_0,\x_0)$.
\ep

The result of \rl{5.1} can be restated as a uniform estimate for the functions
$u(t)\in C^2(\so)$. Since, moreover, the principal curvatures of the flow
hypersurfaces are not only bounded, but also uniformly bounded away from zero,
in view of \re{5.4} and the assumption that $F$ vanishes on $\pa \C_+$, we
conclude that $F$ is uniformly elliptic on $M(t)$.

\section{Convergence to a stationary solution}\las 6

We are now ready to prove \rt{0.2}. Let $M(t)$ be the flow with initial
hypersurface $M_0=M_2$. Let us look at the scalar version of the flow \re{3.5}
\begin{equation}\lae{6.1}
\pde ut=-e^{-\psi}v(\F-\tilde f).
\end{equation}
This is  a scalar parabolic differential equation defined on the cylinder
\begin{equation}
Q_{T^*}=[0,T^*)\times \so
\end{equation}
with initial value $u(0)=u_2\in C^{4,\alpha}(\so)$. In view of the a priori estimates,
which we have established in the preceding sections, we know that
\begin{equation}
{\abs u}_\low{2,0,\so}\le c
\end{equation}
and
\begin{equation}
\F(F)\,\tup{is uniformly elliptic in}\,u
\end{equation}
independent of $t$. Moreover, $\F(F)$ is concave, and thus, we can apply
the regularity results of \ci[Chapter 5.5]{nk} to conclude that uniform
$C^{2,\alpha}$-estimates are valid, leading further to uniform $C^{4,\alpha}$-estimates
due to the regularity results for linear operators.

Therefore, the maximal time interval is unbounded, i.e. $T^*=\un$.

Now, integrating \re{6.1} with respect to $t$, and observing that the right-hand
side is non-positive, yields
\begin{equation}
u(0,x)-u(t,x)=\int_0^te^{-\psi}v(\F-\tilde f)\ge c\int_0^t(\F-\tilde f),
\end{equation}
i.e.,
\begin{equation}
\int_0^\un \abs{\F-\tilde f}<\un\qq\A\msp x\in \so
\end{equation}
Hence, for any $x\in\so$ there is a sequence $t_k\rightarrow \un$ such that
$(\F-\tilde f)\rightarrow 0$.

On the other hand, $u(\cdot,x)$ is monotone decreasing and therefore
\begin{equation}
\lim_{t\rightarrow \un}u(t,x)=\tilde u(x)
\end{equation}
exists and is of class $C^{4,\alpha}(\so)$ in view of the a priori estimates. We, finally,
conclude that $\tilde u$ is a stationary solution of our problem, and that
\begin{equation}
\lim_{t\rightarrow \un}(\F-\tilde f)=0.
\end{equation}

\end{document}